\documentclass{jams-l}

\usepackage{graphicx}
\usepackage{mathtools}
\usepackage{xcolor}

\newtheorem{theorem}{Theorem}[section]
\newtheorem{lemma}[theorem]{Lemma}

\theoremstyle{definition}
\newtheorem{definition}[theorem]{Definition}
\newtheorem{example}[theorem]{Example}

\theoremstyle{remark}

\numberwithin{equation}{section}

\begin{document}

\title[AMS Article Template]{Carleman linearization and systems of arbitrary depth polynomial recursions}

\author{Mikołaj Myszkowski}
\address{Parks Rd, Oxford OX1 3PJ, UK}
\email{mikolaj.myszkowski@hertford.ox.ac.uk}


\date{}

\dedicatory{}


\begin{abstract}
New approach to systems of polynomial recursions is developed based on the Carleman linearization procedure. The article is divided into two main sections: firstly, we focus on the case of uni-variable depth-one polynomial recurrences. Subsequently, the systems of depth-one polynomial recurrence relations are discussed. The corresponding transition matrix is constructed and upper triangularized. Furthermore, the powers of the transition matrix are calculated using the back substitution procedure. The explicit expression for a solution to a broad family of recurrence relations is obtained. We investigate to which recurrences the framework can be applied and construct a sufficient conditions for the method to work. It is shown how introduction of auxiliary variables can be used to reduce arbitrary depth systems to the depth-one system of recurrences dealt with earlier. Finally, the limitations of the method are discussed, outlining possible directions of future research.

\end{abstract}

\maketitle
\section{Introduction}
Recurrence relations arise in various fields of mathematics and prove to be a powerful tool for studying probability \cite{Book_5} and combinatorics \cite{Book_13}. It is often the case that construction of a solution to a mathematical problem simplifies to solving a recurrence relation. The theory of linear recurrences is well studied, leading to various results such as the celebrated Binet's formula or the characteristic polynomial method \cite{Book_1,Book_2}. The dynamics of non linear recurrence sequences is of special interest, as there are various examples of non linear (in particular, polynomial) maps exhibiting chaotic behaviour. The examples include the logistic, quadratic and exponential maps, as well as the Ricatti recurrence \cite{Article_9, Article_2, Article_7, Article_8, Article_1, Article_20}.
\begin{definition}
\cite{Article_6} Let ${u_{i}^{k}\in \mathbb{C}, i\geq 0}$ be ${k}$ sequences of complex numbers such that for every ${i\geq n}$ the following relation holds:
\begin{equation}
\begin{cases}
u_{i}^{1} & =F_{1}\left(u_{i-1}^{1},..., u_{i-1}^{k}, u_{i-2}^{1},..., u_{i-2}^{k},...,u_{i-n}^{1},..., u_{i-n}^{k}\right)\\
u_{i}^{2} & =F_{2}\left(u_{i-1}^{1},..., u_{i-1}^{k}, u_{i-2}^{1},..., u_{i-2}^{k},...,u_{i-n}^{1},..., u_{i-n}^{k}\right)\\
&\vdots\\
u_{i}^{k} & =F_{k}\left(u_{i-1}^{1},..., u_{i-1}^{k}, u_{i-2}^{1},..., u_{i-2}^{k},...,u_{i-n}^{1},..., u_{i-n}^{k}\right)\\
\end{cases}
\end{equation}
Then, we say that sequences ${u_{i}^{k}\in \mathbb{C}}$ are a solution to the depth-n system of recurrences defined by functions ${F_{i}}$.
\end{definition}
In case when functions ${F_{i}}$ are linear, the system of equations (1.1) can be rewritten in the matrix form \cite{Book_6}. If the auxiliary vectors ${{\bf y}^{i},i\geq n}$ and ${{\bf e}_{j}}$ are defined as:
\begin{equation}
{{\bf y}^{i}}=[1,u_{i}^{1},..., u_{i}^{k}, u_{i-1}^{1},..., u_{i-1}^{k}, u_{i-2}^{1},..., u_{i-2}^{k},..., u_{i-n+1}^{1},..., u_{i-n+1}^{k}]^{T}
\end{equation}
\begin{equation}
{{\bf e}_{j}}=\left[\delta_{0,j}, \delta_{1,j},...,\delta_{nk,j}\right]
\end{equation}
then the system of recurrences (1.1) can be rewritten as recurrence relation for the sequence of vectors ${{\bf y}^{i}}$:
\begin{equation}
{{\bf y}^{i}}=M{{\bf y}^{i-1}}
\end{equation}
where ${(1+nk)\times (1+nk)}$ matrix M is called the transition matrix associated with (1.1). By induction, the solution to (1.4) is given by powers of matrix M:
\begin{equation}
{{\bf y}^{i}}=M^{i+1-n}{{\bf y}^{n-1}}
\end{equation}
where ${{\bf y}^{n-1}}$ is defined by initial conditions. Equation (1.5) can be then multiplied by vectors ${{\bf e}_{j}}$ from the left, hence extracting the original sequences from ${{\bf y}^{i}}$:
\begin{equation}
u_{i}^{k}={{\bf e}_{k}}M^{i+1-n}{{\bf y}^{n-1}}
\end{equation}
The above method of solving linear recurrences is fundamentally equivalent to the method of characteristic polynomials \cite{Book_8}.\\\indent
Even if functions ${F_{i}}$ are non linear (but analytic), the solution to the system still can be formulated in terms of powers of appropriately chosen transition matrix, in analogy with the linear counterpart \cite{Article_3}. The method of describing dynamical systems (recursions in particular) in terms of infinite-dimensional linear algebra is referred to as Carleman embedding or Carleman linearization \cite{Book_3, Article_5}.

\begin{theorem}
\cite{Article_5}\cite{Article_1} Let ${k=1}$, and ${F_{1}=F_{1}\left(u_{i-1}^{1}\right)}$ be analytic in ${u_{i-1}^{1}}$. Let the (infinite) auxiliary vectors be defined as:
\begin{equation}
{{\bf y}^{i}}=\left[1, u_{i}^{1}, \left(u_{i}^{1}\right)^{2}, \left(u_{i}^{1}\right)^{3},...\right]^{T}, 
{{\bf e}_{j}}=\left[\delta_{0,j}, \delta_{1,j}, \delta_{2,j},...\right]
\end{equation}
Then, the solution can be written in the form:
\begin{equation}
u_{i}^{1}={{\bf e}_{1}}M^{i}{{\bf y}^{0}}
\end{equation}
where ${{\bf y}^{0}}$ is specified by initial conditions, and the infinite dimensional matrix ${M}$ is given by coefficients of Maclaurin expansion of ${F_{1}}$:

\begin{equation}
M_{ab}=\frac{1}{b!}\frac{d^{b}}{dx^{b}}\left(F_{1}(x)\right)^{a}\biggr\rvert_{x=0}
\end{equation}
\end{theorem}
Theorem 1.2 allows for an analytic function to be iterated by computing powers of the infinite transition matrix associated with the function. In practice, however, it is troublesome to give an explicit expression for the powers of matrix M \cite{Article_11}. Nevertheless, in some cases the transition matrix can be diagonalized and the corresponding eigenvectors can be found \cite{Article_5}. We focus on systems of recurrence relations (1.1) with ${F_{i}}$ restricted to be finite degree polynomials. The article is divided into two main sections.\\\indent
Firstly, the k=1 (uni-variable) depth-one case is discussed, leading to the exact expression for the solution in terms of the initial condition and coefficients of polynomials ${F_{i}}$. As it turns out, in most cases the transition matrix can be diagonalized by a suitable choice of variables.\\\indent
The method of the Carleman embedding is then applied to the multi-variable depth one recurrence. With minor modifications of the auxiliary vectors, the transition matrix can be compactly written as a block matrix. Again, M is diagonalized and the eigenvectors are found. The main theorem of the article is given, solving the recursion.\\\indent
Finally, we discuss the general case (1.1). Provided that ${F_{i}}$ are polynomial in ${u_{i}^{k}}$, the system of equations (1.1) is embedded in an infinite linear space. It is shown that by introducing appropriate auxiliary variables, any multi-variable arbitrary depth recurrence can be reduced to multi-variable depth-one recurrence. Therefore, the recurrence (1.1) can be solved using the main theorem presented in Section 3.\\\indent
Throughout the article we assume that the functions ${F_{i}}$ are polynomial in all variables. For convenience, the first elements of vectors are indexed from ${0}$ instead of ${1}$. A similar rule applies to matrices.

\section{uni-variable recurrences of depth one}
In the uni-variable case, there is only one sequence ${u_{i}^{1}\equiv u_{i}}$ and one function ${F_{i}\equiv F}$ (we drop the unnecessary indices for simplicity). Furthermore, ${F}$ is a function of a single variable. The recurrence relations of this kind can be written as:
\begin{equation}
u_{i}=F\left(u_{i-1}\right)=\sum_{j=0}^{m}c_{j}\left(u_{i-1}\right)^{j}
\end{equation}
where ${c_{j}}$ are constant coefficients of the polynomial ${F}$. The recurrence (2.1) is of special interest, as it includes many well-known nonlinear chaotic maps, i.e. the logistic map, that is used as an example at the end of this section \cite{Article_19}.\\\indent
Using Theorem 1.2 and the multinomial expansion theorem \cite{Article_12}, the infinite dimensional transition matrix can be obtained:
\begin{equation}
\begin{aligned}
T_{ab}&=\frac{1}{b!}\frac{d^{b}}{dx^{b}}\left(F_{1}(x)\right)^{a}\biggr\rvert_{x=0}=\frac{1}{b!}\frac{d^{b}}{dx^{b}}\left(\sum_{j=0}^{m}c_{j}x^{j}\right)^{a}\biggr\rvert_{x=0}=\\
&=\frac{1}{b!}\frac{d^{b}}{dx^{b}}\sum_{\sum_{l=0}^{m}k_{l}=a}\binom{a}{k_{0},k_{1},...,k_{m}}\prod_{t=0}^{m}\left(c_{t}x^{t}\right)^{k_{t}}\biggr\rvert_{x=0}
\end{aligned}
\end{equation}
where ${k_{l}}$ are assumed to be non-negative integers. The condition ${\sum_{l=0}^{m}k_{l}=a}$ denotes the fact that the bigger sum runs only over such non-negative integers ${k_{l}}$ that sum up to ${a}$. This further simplifies to:
\begin{equation}
\begin{aligned}
T_{ab}&=\frac{1}{b!}\frac{d^{b}}{dx^{b}}\sum_{\sum_{l=0}^{m}k_{l}=a}\binom{a}{k_{0},k_{1},...,k_{m}}\prod_{t=0}^{m}\left(c_{t}x^{t}\right)^{k_{t}}\biggr\rvert_{x=0}=\\
&=\sum_{\sum_{l=0}^{m}k_{l}=a,\ \sum_{l=0}^{m}lk_{l}=b}\binom{a}{k_{0},k_{1},...,k_{m}}\prod_{t=0}^{m}c_{t}^{k_{t}}=\\
&=\sum_{\sum_{l=0}^{m}k_{l}=a,\ \sum_{l=0}^{m}lk_{l}=b}\frac{a!}{\prod_{s=0}^{m}k_{s}!}\prod_{t=0}^{m}c_{t}^{k_{t}}
\end{aligned}
\end{equation}
Again, the conditions under the bigger sum narrow down the possible combinations of ${k_{l}}$ to those that obey ${\sum_{l=0}^{m}k_{l}=a,\ \sum_{l=0}^{m}lk_{l}=b}$. If no such integers exist, the sum vanishes and ${T_{ab}=0}$. The transition matrix (2.3), together with the auxiliary vectors from equation (1.7) define the Carleman embedding of the recurrence (2.1) in the infinite dimensional linear space. We now prove several facts that will become useful later on in calculating powers of ${T_{ab}}$.
\begin{theorem}
Let ${T_{ab}}$ be the infinite dimensional transition matrix corresponding to the recurrence (2.1). Then, the following holds true:
\begin{enumerate}
  \item ${T_{0b}=\delta_{0,b}}$
  \item If ${c_{0}=0}$, ${T_{ab}}$ is upper triangular
\end{enumerate}
\end{theorem}

\begin{proof}
${\bf (1)}$ We obtain:
\begin{equation}
T_{0b}=\sum_{\sum_{l=0}^{m}k_{l}=0,\ \sum_{l=0}^{m}lk_{l}=b}\binom{0}{k_{0},k_{1},...,k_{m}}\prod_{t=0}^{m}c_{t}^{k_{t}}
\end{equation}
However, since ${k_{l}\geq 0}$ and ${\sum_{l=0}^{m}k_{l}=0}$, all of the integers ${k_{l}}$ have to vanish. Therefore ${\sum_{l=0}^{m}lk_{l}=b=0}$. If ${b\neq 0}$, there are no ${k_{l}}$ satisfying the conditions under the sum, and ${T_{0b}=0}$. If ${b=0}$:
\begin{equation}
   T_{00}=\binom{0}{0,...,0}\prod_{t=0}^{m}c_{t}^{0}=1
\end{equation}
{\bf (2)} If ${c_{0}=0}$, then the lowest non-zero power of ${x}$ appearing in polynomial ${F(x)}$ is ${x^{1}}$. By extension, the lowest power of ${x}$ appearing in powers ${F^{a}}$ is of order ${a}$. As a result:
\begin{equation}
T_{ab}=\frac{1}{b!}\frac{d^{b}}{dx^{b}}\left(\sum_{j=1}^{m}c_{j}x^{j}\right)^{a}\biggr\rvert_{x=0}=\frac{1}{b!}\frac{d^{b}}{dx^{b}}\left(c_{1}^{a}x^{a}+...\right)\biggr\rvert_{x=0}=0
\end{equation}
for ${b<a}$.
\end{proof}
It is worth noting that the above mentioned polynomials have also been studied in the setting of graphs \cite{Article_21}. The fact that the transition matrix can take the upper triangular form makes diagonalization feasible. One is always free to work with new shifted variables ${u_{i}'=u_{i}-d}$, which satisfy a new recurrence:
\begin{equation}
u_{i}'=F'\left(u_{i-1}'\right)=\sum_{j=0}^{m}c_{j}'\left(u_{i-1}'\right)^{j}
\end{equation}
The coefficients ${c_{j}'}$ can be obtained by substituting ${u_{i}'=u_{i}-d}$ into (2.1):
\begin{equation}
\begin{aligned}
u_{i}'+d&=\sum_{j=0}^{m}c_{j}\left(u_{i-1}'+d\right)^{j}=\sum_{j=0}^{m}c_{j}\sum_{l=0}^{j}\binom{j}{l}d^{l}u_{i-1}'^{j-l}=\\
&=\sum_{l=0}^{m}\sum_{j=l}^{m}c_{j}\binom{j}{l}d^{j-l}u_{i-1}'^{l}=\sum_{j=0}^{m}c_{j}'u_{i-1}'^{j}+d
\end{aligned}
\end{equation}
Which implies:
\begin{equation}
c_{0}'=\sum_{j=0}^{m}c_{j}d^{j}-d, \ c_{l}'=\sum_{j=l}^{m}c_{j}\binom{j}{l}d^{j-l} \ for \ l\neq 0
\end{equation}
Therefore, new coefficients ${c_{j}'}$ are functions of old coefficients ${c_{j}}$ and the shift parameter ${d}$. It turns out that there always exists ${d}$ such that the transition matrix of the transformed recurrence becomes upper triangular.
\begin{lemma}
Let the recurrence on ${u_{i}}$ be defined as in (2.1). For every ${F(u_{i-1})}$ there exists ${d\in \mathbb{C}}$ such that the transition matrix ${T_{ab}'}$ corresponding to the new sequence ${u_{i}'=u_{i}-d}$ is upper triangular.
\end{lemma}
\begin{proof}
The new sequence ${u_{i}'}$ is obtained by performing a shift transformation with a parameter ${d}$. As already derived in (2.7), the recurrence of the shifted sequence is:
\begin{equation}
u_{i}'=\sum_{j=0}^{m}c_{j}'u_{i-1}'^{j}
\end{equation}
where the constant term (term near 0-th power of ${u_{i-1}'}$) reads:
\begin{equation}
c_{0}'=\sum_{j=0}^{m}c_{j}d^{j}-d
\end{equation}
Theorem 2.1 states that the transition matrix ${T_{ab}'}$ corresponding to the shifted recurrence is upper triangular if the constant term ${c_{0}'=0}$ vanishes. This leads to the following polynomial equation in ${d}$:
\begin{equation}
0=\sum_{j=0}^{m}c_{j}d^{j}-d
\end{equation}
By the fundamental theorem of algebra, there always exists ${d\in C}$ such that ${c_{0}'=0}$ \cite{Article_16}. Therefore, the transition matrix ${T_{ab}}$ can be always brought to an upper triangular form by a suitable shift of variables.
\end{proof}
As a result, we can consider only the recurrences satisfying ${c_{0}=0}$ without loss of generality. \\\indent
Because the transition matrix is upper triangular, the diagonal elements ${T_{aa}}$ are it's eigenvalues \cite{Article_17}. Using (2.2):
\begin{equation}
\lambda_{a}=T_{aa}=\frac{1}{a!}\frac{d^{a}}{dx^{a}}\left(\sum_{j=1}^{m}c_{j}x^{j}\right)^{a}\biggr\rvert_{x=0}=\frac{1}{a!}\frac{d^{a}}{dx^{a}}\left(c_{1}^{a}x^{a}+...\right)\biggr\rvert_{x=0}=c_{1}^{a}
\end{equation}
The transition matrix is diagonalizable if all its eigenvalues are different, i.e. when ${c_{1}^{a}=c_{1}^{b}\implies a=b}$ \cite{Book_10}. Since in general ${c_{1}\in \mathbb{C}}$, this condition can be rewritten as ${c_{1}\neq 0\wedge(|c_{1}|\neq 1\vee arg(c_{1})\notin \mathbb{Q})}$. It should be noted that the right hand side of equation (2.11) may have more than one distinct root, which means that there exist recurrences that can be shift transformed from ${c_{1}= 0\vee(|c_{1}|=1\wedge arg(c_{1})\in \mathbb{Q})}$ to ${c_{1}\neq 0\wedge(|c_{1}|\neq 1\vee arg(c_{1})\notin \mathbb{Q})}$.
\begin{example}
Let ${F(x)=x^{3}+2x^{2}+x}$. The recurrence reads:
\begin{equation}
u_{i}=F(u_{i-1})=(u_{i-1})^{3}+2(u_{i-1})^{2}+u_{i-1}
\end{equation}
The linear coefficient ${c_{1}=1}$, hence there is no guarantee that the corresponding transition matrix is diagonalizable. Nevertheless, by performing the shift transformation with ${d=-2}$, we obtain (using formulas (2.9)):
\begin{equation}
u_{i}'=(u_{i-1}')^{3}-4(u_{i-1}')^{2}+5u_{i-1}'
\end{equation}
with ${c'_{0}=0, |c'_{1}|\neq 1\neq 0}$.
\end{example}
The above example shows that the condition ${c_{1}\neq 0\wedge(|c_{1}|\neq 1\vee arg(c_{1})\notin \mathbb{Q})}$ is sufficient, but not necessary for the diagonalization of $T_{ab}$. Some transition matrices that are hard to deal with can be diagonalized after shift of variables, thus the recurrence can be solved.\\\indent
The next two Lemmas regard the procedure of diagonalization of infinite upper triangular matrices. The set of linearly independent eigenvectors is obtained. Even though the matrix is infinite, the eigenvector corresponding to a-th eigenvalue ${\lambda_{a}}$ has at most ${a+1}$ non-zero components. An explicit expression for the inverse of an upper triangular matrix is given (provided that it exists), which allows for calculation of the inverse of the corresponding modal matrix. The transition matrix is then diagonalized, and the powers of the matrix are subsequently calculated.
\begin{lemma}
Let ${M_{ab}}$ be an infinite upper triangular matrix with ${\forall a \ \lambda_{a}=M_{aa}\neq 0}$ and different eigenvalues, i.e. ${\lambda_{a}=\lambda_{b}\implies a=b}$. Let ${v_{a}^{b}}$ denote the b-th component of the (not necessary normalized) eigenvector corresponding to the a-th eigenvalue ${\lambda_{a}}$. Then:
\begin{enumerate}
  \item ${\forall b>a \ v_{a}^{b}=0}$
  \item ${\forall b=a \ v_{a}^{b}=1}$
  \item ${\forall l_{0}<l_{p+1}}$:
\begin{equation}
v_{l_{p+1}}^{l_{0}}=\sum_{p=0}^{l_{p+1}-l_{0}-1} \ \sum_{l_{0}<l_{1}<...<l_{p}<l_{p+1}}(-1)^{p+1}\prod_{j=0}^{p}\frac{M_{l_{j},l_{j+1}}}{M_{l_{j},l_{j}}-\lambda_{l_{p+1}}}
\end{equation}
\end{enumerate}
\end{lemma}
\begin{proof}
The eigenvector equation ${Mv_{a}=\lambda_{a}v_{a}}$ in the matrix form reads:
\begin{equation}
\begingroup
\setlength\arraycolsep{3pt}
\begin{bmatrix}
M_{00}-\lambda_{a} & M_{01} & M_{02} & M_{03} & M_{04} & M_{05} & \dots\\
0 & \ddots & \vdots & \vdots & \vdots & \vdots & \\
0 & \dots & M_{a-1,a-1}-\lambda_{a} & M_{a-1,a} & M_{a-1,a+1} & M_{a-1,a+2} & \dots \\
0 & \dots & 0 & 0 & M_{a,a+1} & M_{a,a+2} & \dots \\
0 & \dots & 0 & 0 & M_{a+1,a+1} & M_{a+1,a+2} & \dots \\
\vdots &  & \vdots & \vdots & \vdots & \vdots & \ddots \\
\end{bmatrix}
\endgroup
\begin{bmatrix}
v_{a}^{0}\\
\vdots\\
v_{a}^{a-1}\\
v_{a}^{a}\\
v_{a}^{a+1}\\
\vdots \\
\end{bmatrix}
=0
\end{equation}
Lets consider a finite upper triangular system of equations that arises from truncating the above matrix at ${(a+1)\times (a+1)}$:
\begin{equation}
\begin{bmatrix}
M_{00}-\lambda_{a} & M_{01} & \dots & M_{0a}\\
0 & \ddots & \vdots &\vdots & \\
0 & \dots & M_{a-1,a-1}-\lambda_{a} & M_{a-1,a} & \\
0 & \dots & 0 & 0 \\
\end{bmatrix}
\begin{bmatrix}
v_{a}^{0} \\
\vdots \\
v_{a}^{a-1} \\
v_{a}^{a} \\
\end{bmatrix}
=0
\end{equation}
Since all of the eigenvalues are different, the diagonal elements ${\forall b\neq a \ \lambda_{b}-\lambda_{a}\neq 0}$ and the matrix ${M}$ is diagonalizable. The determinant of (2.16) vanishes, and thus the solution is not unique (up to a scaling factor), as expected. The explicit expression for ${v_{a}^{b}}$ as a function of ${v_{a}^{a}}$ can be obtained by the back substitution method \cite{Book_11}:
\begin{equation}
v_{l_{p+1}}^{l_{0}}=v_{l_{p+1}}^{l_{p+1}}\sum_{p=0}^{l_{p+1}-l_{0}-1} \ \sum_{l_{0}<l_{1}<...<l_{p}<l_{p+1}}(-1)^{p+1}\prod_{j=0}^{p}\frac{M_{l_{j},l_{j+1}}}{M_{l_{j},l_{j}}-\lambda_{l_{p+1}}}
\end{equation}
For simplicity, we set the scaling factor ${v_{l_{p+1}}^{l_{p+1}}=1}$ for all ${l_{p+1}}$. It can be seen that ${l_{0}>l_{p+1}}$ would contradict with ${l_{0}<l_{1}<...<l_{p}<l_{p+1}}$, hence (2.19) vanishes and the identity (1) holds (for the eigenvectors of the truncated matrix). Since we know that ${\forall b>a \ v_{a}^{b}=0}$, the eigenvectors of the truncated matrix (2.18) with zeros substituted for the rest of the components are also the solution to the infinite equation (2.17). Therefore, expressions (1), (2), (3) and (2.16) also hold for the matrix ${M}$.
\end{proof}
The following formula is an explicit expression for the inverse of an upper triangular matrix:
\begin{lemma}
Let ${U}$ be an infinite upper triangular matrix with complex entries:
\begin{equation}
U_{km}=
\begin{bmatrix}
U_{00} & U_{01} & U_{02} & U_{03} & \dots \\
0 & U_{11} & U_{12} & U_{13} & \dots \\
0 & 0 & U_{22} & U_{23} & \dots \\
0 & 0 & 0 & U_{33} & \dots \\
\vdots & \vdots & \vdots & \vdots & \ddots \\
\end{bmatrix}
\end{equation}
and let ${\forall n \ U_{nn}\neq 0}$. The explicit expression for the inverse of ${U}$ is given by:
\begin{equation}
U_{km}^{-1}=
\begin{cases}
0 \ for \ k>m \\
\frac{1}{U_{km}} \ for \ k=m \\
\end{cases}
\end{equation}
and:
\begin{equation}
U_{l_{0},l_{p+1}}^{-1}=\frac{1}{U_{l_{p+1},l_{p+1}}}\sum_{p=0}^{l_{p+1}-l_{0}-1} \ \sum_{l_{0}<l_{1}<...<l_{p}<l_{p+1}}(-1)^{p+1}\prod_{j=0}^{p}\frac{U_{l_{j},l_{j+1}}}{U_{l_{j},l_{j}}}
\end{equation}
for ${l_{0}<l_{p+1}}$.
\end{lemma}
\begin{proof}
The inverse of the upper triangular matrix is also upper triangular \cite{Article_13}, therefore the equation ${UU^{-1}=I}$ reads:
\begin{equation}
\begin{bmatrix}
U_{00} & U_{01} & U_{02} & U_{03} & \dots \\
0 & U_{11} & U_{12} & U_{13} & \dots \\
0 & 0 & U_{22} & U_{23} & \dots \\
0 & 0 & 0 & U_{33} & \dots \\
\vdots & \vdots & \vdots & \vdots & \ddots \\
\end{bmatrix}
\begin{bmatrix}
U_{00}^{-1} & U_{01}^{-1} & U_{02}^{-1} & U_{03}^{-1} & \dots \\
0 & U_{11}^{-1} & U_{12}^{-1} & U_{13}^{-1} & \dots \\
0 & 0 & U_{22}^{-1} & U_{23}^{-1} & \dots \\
0 & 0 & 0 & U_{33}^{-1} & \dots \\
\vdots & \vdots & \vdots & \vdots & \ddots \\
\end{bmatrix}
=I
\end{equation}
For the diagonal elements, we have:
\begin{equation}
U_{nn}^{-1}=\frac{1}{U_{nn}}
\end{equation}
Elements further off the diagonal can be obtained by the back substitution algorithm. Using (2.23):
\begin{equation}
U_{nn}U_{n,n+1}^{-1}+U_{n,n+1}U_{n+1,n+1}^{-1}=0
\end{equation}
Which together with (2.24) gives:
\begin{equation}
U_{n,n+1}^{-1}=-\frac{U_{n,n+1}U_{n+1,n+1}^{-1}}{U_{nn}}=-\frac{U_{n,n+1}}{U_{nn}U_{n+1,n+1}}
\end{equation}
By repeating the same process for the rest of the elements of ${U_{km}^{-1}}$, the relation (2.22) for ${k<m}$ is obtained.
\end{proof}
The powers of the transition matrix can be calculated by diagonalizing ${T}$:
\begin{equation}
T^{n}=\left(PDP^{-1}\right)^{n}=PD^{n}P^{-1}
\end{equation}
where ${D}$ denotes the infinite matrix with the eigenvalues of ${T}$ as diagonal entries, while P is the infinite modal matrix composed of eigenvectors of ${T}$:
\begin{equation}
D=diag[\lambda_{0}, \lambda_{1}, \lambda_{2},...], \ P =[ {\bf v}_{0}, {\bf v}_{1}, {\bf v}_{2},...]
\end{equation}
We now prove the main theorem of this section:
\begin{theorem}
Let ${u_{i}\in \mathbb{C}, i\geq 0}$ be a sequence of numbers satisfying the recurrence relation (2.1). Let ${F(u_{i-1})=\sum_{j=0}^{m}c_{j}\left(u_{i-1}\right)^{j}}$ be polynomial in ${u_{i-1}}$ such that the corresponding polynomial equation:
\begin{equation}
\sum_{j=0}^{m}c_{j}d^{j}=d
\end{equation}
has a root ${d\in \mathbb{C}}$ satisfying:
\begin{equation}
\left\lvert \sum_{j=1}^{m}jc_{j}d^{j-1} \right\rvert\neq 0\wedge \left(\left\lvert \sum_{j=1}^{m}jc_{j}d^{j-1} \right\rvert\neq 1\vee arg\left(\sum_{j=1}^{m}jc_{j}d^{j-1}\right)\notin \mathbb{Q}\right)
\end{equation}
Then, the solution to the recurrence (2.1) is given by:
\begin{equation}
u_{i}=d+\sum_{l=1}^{\infty}\sum_{j=1}^{l}v'^{1}_{j}c'^{ij}_{1}P'^{-1}_{jl}(u_{0}-d)^{l}
\end{equation}
where:
\begin{equation}
\begin{aligned}
&P'^{-1}_{l_{0},l_{p+1}}=
\begin{cases}
0 \ for \ l_{0}>l_{p+1} \\
1 \ for \ l_{0}=l_{p+1} \\
\sum_{p=0}^{l_{p+1}-l_{0}-1} \ \sum_{l_{0}<l_{1}<...<l_{p}<l_{p+1}}(-1)^{p+1}\prod_{j=0}^{p}\frac{v'^{l_{j}}_{l_{j+1}}}{v'^{l_{j}}_{l_{j}}}\  
for \ l_{0}<l_{p+1}
\end{cases}
\\
&v'^{l_{0}}_{l_{p+1}}=
\begin{cases}
0 \ for \ l_{0}>l_{p+1} \\
1 \ for \ l_{0}=l_{p+1} \\
\sum_{p=0}^{l_{p+1}-l_{0}-1} \ \sum_{l_{0}<l_{1}<...<l_{p}<l_{p+1}}(-1)^{p+1}\prod_{j=0}^{p}\frac{T'_{l_{j},l_{j+1}}}{T'_{l_{j},l_{j}}-\lambda_{l_{p+1}}} \ for \ l_{0}<l_{p+1}

\end{cases}
\\
&T_{ab}'=\sum_{\sum_{l=0}^{m}k_{l}=a,\ \sum_{l=0}^{m}lk_{l}=b}\frac{a!}{\prod_{s=0}^{m}k_{s}!}\prod_{t=0}^{m}c_{t}'^{k_{t}}
\\
&c_{0}'=0, \ c_{l}'=\sum_{j=l}^{m}c_{j}\binom{j}{l}d^{j-l} \ for \ l\neq 0
\end{aligned}
\end{equation}
\end{theorem}
\begin{proof}
Because the equation (2.29) is equivalent to requiring ${c_{0}'=0}$ (${c_{0}'}$ is the free term of the polynomial ${F'}$ obtained after shift of variables by ${d}$), the variables ${u_{i}'=u_{i}-d}$ are the ones in which the free term vanishes. In the same way, from (2.9) we have:
\begin{equation}
c_{1}'=\sum_{j=1}^{m}c_{j}\binom{j}{1}d^{j-1}=\sum_{j=1}^{m}jc_{j}d^{j-1}
\end{equation}
Therefore the condition (2.30) expresses the fact that all of the powers of ${c_{1}'}$ have different values.\\\indent
The transition matrix ${T}$ is (2.3):
\begin{equation}
T_{ab}'=\sum_{\sum_{l=0}^{m}k_{l}=a,\ \sum_{l=0}^{m}lk_{l}=b}\frac{a!}{\prod_{s=0}^{m}k_{s}!}\prod_{t=0}^{m}c_{t}'^{k_{t}}
\end{equation}
By Lemma 2.2, ${c_{0}'=0}$ and the infinite transition matrix is upper triangular, while (2.30) implies that ${T}$ is diagonalizable. Furthermore, by Lemma 2.4 the a-th eigenvector has at most ${a+1}$ non-zero elements, and the (infinite) matrix ${P}$ can be chosen to be upper triangular (the inverse of ${P}$ is then upper triangular as well, Lemma 2.5). The powers of transition matrix ${T}$ can be obtained from (2.27). Hence, using (1.8), (2.34) and Lemma 2.4:
\begin{equation}
\begin{aligned}
&u_{i}'={\bf e_{1}}T'^{i}{\bf y}'^{0}={\bf e_{1}}P{D'^{i}}P^{-1}{\bf y}'^{0}=
\begin{bmatrix}
0 \\
1 \\
0 \\
\vdots \\
\end{bmatrix}
^{T}
\begin{bmatrix}
v'^{0}_{0} & v'^{0}_{1} & v'^{0}_{2} & \dots \\
0 & v'^{1}_{1} & v'^{1}_{2} & \dots \\
0 & 0 &v'^{2}_{2} & \dots \\
\vdots & \vdots & \vdots & \ddots \\
\end{bmatrix}
\times \\
&\times
\begin{bmatrix}
1 & 0 & 0 & \dots \\
0 & c'_{1} & 0 & \dots \\
0 & 0 & c'^{2}_{1} & \dots \\
\vdots & \vdots & \vdots & \ddots \\
\end{bmatrix}
^{i}
\begin{bmatrix}
P'^{-1}_{00} & P'^{-1}_{01} & P'^{-1}_{02} & \dots \\
0 & P'^{-1}_{11} &  P'^{-1}_{12} & \dots \\
0 & 0 & P'^{-1}_{22} & \dots \\
\vdots & \vdots & \vdots & \ddots \\
\end{bmatrix}
\begin{bmatrix}
1 \\
u_{0}' \\
u'^{2}_{0} \\
\vdots \\
\end{bmatrix}
\end{aligned}
\end{equation}
Where ${P_{ab}^{-1}}$ denote the elements of the inverse of ${P}$ that can be calculated using Lemma 2.5, while the formula for ${v'^{j}_{l}}$ is given in (2.16). The ${i}$-th power of the diagonal matrix is just a matrix with ${i}$-th power of its elements on the diagonal. Carrying out the matrix multiplication, we obtain:
\begin{equation}
u_{i}'=\sum_{l=1}^{\infty}u'^{l}_{0}\sum_{j=1}^{l}v'^{1}_{j}\left(c'^{j}_{1}\right)^{i}P'^{-1}_{jl}=\sum_{l=1}^{\infty}\sum_{j=1}^{l}v'^{1}_{j}c'^{ij}_{1}P'^{-1}_{jl}u'^{l}_{0}
\end{equation}
By retrieving the original sequence ${u_{i}=u_{i}'+d}$ and substituting back the original coefficients ${c_{j}}$ for the shifted parameters ${c_{j}'}$, we arrive at the result.
\end{proof}
Theorem 2.6 allows for solving a large family of nonlinear recurrences. It is apparent from (2.31) that the general solution has a qualitative form:
\begin{equation}
u_{i}=\sum_{j=0}^{\infty}f_{j}(i)u_{0}^{j}
\end{equation}
where ${f_{j}}$ are some (in general complex) functions of number of iterations ${i}$. Therefore, the Theorem 2.6 can be used to determine the form of functions ${f_{j}}$. As an example, the results are applied to the well-known logistic mapping:

\begin{example}
The logistic map is a recurrence relation of the form \cite{Book_7}:
\begin{equation}
u_{i}=ru_{i-1}-ru_{i-1}^{2}
\end{equation}
with the coefficient ${r\in \mathbb{R}}$. The conditions (2.29) and (2.30) are trivially fulfilled for ${r\neq 1}$ and ${d=0}$. We obtain:
\begin{equation}
P=
\begin{bmatrix}
1 & 0 & 0 & 0 & \dots \\
0 & 1 & \frac{1}{1-r} & \frac{2}{r^{3}-r^{2}-r+1} & \dots \\
0 & 0 & 1 & \frac{2}{1-r} & \dots \\
0 & 0 & 0 & 1 & \dots \\
\vdots & \vdots & \vdots & \vdots & \ddots \\
\end{bmatrix}
, P^{-1}=
\begin{bmatrix}
1 & 0 & 0 & 0 & \dots \\
0 & 1 & \frac{1}{r-1} & \frac{2r}{r^{3}-r^{2}-r+1} & \dots \\
0 & 0 & 1 & \frac{2}{r-1} & \dots \\
0 & 0 & 0 & 1 & \dots \\
\vdots & \vdots & \vdots & \vdots & \ddots \\
\end{bmatrix}
\end{equation}
Thus, by Theorem 2.6 the solution takes the form:
\begin{equation}
u_{i}=r^{i}u_{0}+\frac{r^{i}-r^{2i}}{r-1}u_{0}^{2}+\frac{2rr^{i}-2(r+1)r^{2i}+2r^{3i}}{r^{3}-r^{2}-r+1}u_{0}^{3}+...
\end{equation}
\end{example}
It can be seen that the "coefficient functions" ${f_{j}}$ in front of the powers of initial value ${u_{0}}$ become progressively more complicated as the order of ${u_{0}^{j}}$ increases.
\section{Multi-variable recurrences of depth one}
The treatment of polynomial multi-variable recurrences of depth one follows the same structure as Section 2. With appropriately chosen auxiliary vectors, the recurrence can be Carleman-linearized and represented by an appropriate transition matrix. Let:
\begin{equation}
{\bf z}^{i}=[u_{i}^{1}, u_{i}^{2},...,u_{i}^{k}]^{T}
\end{equation}
be a vector composed of ${k}$ variables ${u_{i}^{k}}$. Then, a general polynomial recurrence of depth one in ${k}$ variables can be compactly written in the form:
\begin{equation}
\begin{aligned}
&u_{i}^{p+1}=z_{p}^{i}=F_{p+1}\left(u_{i-1}^{1}, u_{i-1}^{2},...,u_{i-1}^{k}\right)=\prescript{0}{p}C+\sum_{j=1}^{m}\sum_{l_{1},...,l_{j}=0}^{k-1}\prescript{j}{p}C_{l_{1},...,l_{j}}\prod_{s=1}^{j}z_{l_{s}}^{i-1}=\\
&=\prescript{0}{p}C+\sum_{l_{1}=0}^{k-1}\prescript{1}{p}C_{l_{1}}z_{l_{1}}^{i-1}+\sum_{l_{1}=0}^{k-1}\sum_{l_{2}=0}^{k-1}\prescript{2}{p}C_{l_{1},l_{2}}z_{l_{1}}^{i-1}z_{l_{2}}^{i-1}+...
\end{aligned}
\end{equation}
where symbols ${\prescript{j}{p}C_{l_{1},...,l_{j}}}$ represent j+1 dimensional arrays. It is apparent that the description in terms of ${\prescript{j}{p}C_{l_{1},...,l_{j}}}$ is somewhat redundant: every permutation of lower right indices ${l_{1},...,l_{j}}$ is the coefficient of the same combination of powers of variables ${u_{i}^{k}}$. For convenience, we choose ${\prescript{j}{p}C_{l_{1},...,l_{j}}}$ to vanish unless the lower right indices are in the non-decreasing order, i.e. ${\prescript{j}{p}C_{l_{1},...,l_{j}}=0 \ if \ \neg (l_{1}\leq ... \leq l_{j})}$. As a result, the symbol ${\prescript{j}{p}C_{l_{1},...,l_{j}}}$ has only ${k\binom{k+j-1}{j}}$ independent components, when all of the lower indices are taken into account \cite{Article_18}.\\\indent
\begin{theorem}
Let the polynomial multi-variable depth one recurrence be defined as in (3.2), and let the Kronecker powers of ${{\bf z}^{i}}$ be defined as in \cite{Book_8}:
\begin{equation}
\left({\bf z}^{i}\right)^{\alpha}=[u_{i}^{1}\left({\bf z}^{i}\right)^{\alpha-1},u_{i}^{2}\left({\bf z}^{i}\right)^{\alpha-1},...,u_{i}^{k}\left({\bf z}^{i}\right)^{\alpha-1}]^{T}, \ \left({\bf z}^{i}\right)^{0}:=1
\end{equation} 
Then, the Carleman linearization of the system can be defined on infinite auxiliary vectors:
\begin{equation}
{\bf y}^{i}=[\left({\bf z}^{i}\right)^{0},\left({\bf z}^{i}\right)^{1},\left({\bf z}^{i}\right)^{2},...]^{T}
\end{equation}
and the solution to recurrence (3.2) can be written as:
\begin{equation}
u_{i}^{k}={{\bf e}_{k}}T^{i}{{\bf y}^{0}}
\end{equation}
Where ${{\bf y}^{0}}$ is defined by initial conditions and the infinite dimensional transition matrix is:
\begin{equation}
T_{ab}=
\begin{bmatrix}
1 & 0 & 0 & \dots\\
\prescript{0}{p}C & \prescript{1}{p}C_{l_{1}} & \prescript{2}{p}C_{1,l_{1}} & \dots\\
\prescript{0}{p}C \otimes \prescript{0}{p}C & \prescript{0}{p}C \otimes \prescript{1}{p}C_{l_{1}}+\prescript{1}{p}C_{l_{1}} \otimes \prescript{0}{p}C & \dots & \dots\\
\prescript{0}{p}C \otimes \prescript{0}{p}C \otimes \prescript{0}{p}C & \vdots & \ddots & \ddots\\
\vdots & \vdots & \ddots & \ddots\\
\end{bmatrix}
\end{equation}
The Kronecker product in (3.6) acts on the subscripts only.
\end{theorem}
\begin{proof}
Since ${{\bf y}^{i}}$ is composed of consecutive Kronecker powers of ${{\bf z}^{i}}$, any polynomial ${F(u_{i}^{1}, u_{i}^{2},...,u_{i}^{k})}$ in ${k}$ variables can be rewritten as linear combination of elements of ${{\bf y}^{i}}$. The recurrence is defined solely by the arrays ${\prescript{j}{p}C_{l_{1},...,l_{j}}}$, therefore the transition matrix consists of the numbers appearing in ${\prescript{j}{p}C_{l_{1},...,l_{j}}}$. The elements of vector ${{\bf y}^{i}}$ can be identified with the product of elements of ${{\bf z}^{i}}$. From the definitions (3.3), (3.4) and the properties of the Kronecker product we have:
\begin{equation}
\begin{aligned}
y_{l}^{i}=
\begin{cases}
1 \ for \ l=0 \\
\prod_{s=1}^{\lfloor \log_{k}(1+l(k-1)) \rfloor}u_{i}^{(\lfloor (l-\frac{1-k^{s}}{1-k})k^{1-s}\rfloor\%k)+1}
\end{cases}
\\
\prod_{j=1}^{s}u_{i}^{l_{j}}=y_{\frac{1-k^{s}}{1-k}+\sum_{j=1}^{s}k^{s-j}(l_{j}-1)}^{i}=y_{\sum_{j=1}^{s}k^{s-j}l_{j}}^{i}
\end{aligned}
\end{equation} 
The product of elements of ${{\bf z}^{i}}$ can be expressed in terms of elements of ${{\bf z}^{i-1}}$:
\begin{equation}
\prod_{j=0}^{k-1}\left(z_{j}^{i}\right)^{t_{j}}=\prod_{j=0}^{k-1}\left(F_{j}({\bf z}^{i-1})\right)^{t_{j}}
\end{equation}
Analogously to (1.9), the infinite transition matrix can be defined as:
\begin{equation}
\begin{aligned}
T_{ab}=&
R_{b}\left(\prod_{s=1}^{\lfloor \log_{k}(1+b(k-1)) \rfloor}\frac{\partial}{\partial u_{i-1}^{(\lfloor (b-\frac{1-k^{s}}{1-k})k^{1-s}\rfloor\%k)+1}}\right)\cdot\\
\cdot&\left(\prod_{s=1}^{\lfloor \log_{k}(1+a(k-1)) \rfloor}F_{(\lfloor (a-\frac{1-k^{s}}{1-k})k^{1-s}\rfloor\%k)+1}(u_{i-1}^{1},...,u_{i-1}^{k})\right)\biggr\rvert_{u_{i-1}^{k}=0}
\end{aligned}
\end{equation}
The constants ${R_{b}}$ do not depend on ${F_{i}}$, and are defined as a product:
\begin{equation}
R_{b}=\frac{1}{c_{b}}\prod_{s=1}^{k}\frac{1}{a_{s}!}
\end{equation}
where ${a_{s}}$ is the number of times ${u_{i-1}^{s}}$ appears in ${y_{b}^{i-1}}$, while ${c_{b}}$ is the number of times the element ${y^{i}_{b}}$ appears in ${{\bf y}^{i}}$. Equivalently, ${T_{ab}}$ can be written in the matrix form:
\begin{equation}
\begingroup 
\setlength\arraycolsep{1.5pt}
\left[ \begin{array}{@{}*{20}{c}@{}}
1 & 0 & \dots & 0 & \dots \\
\prescript{0}{0}C & \prescript{1}{0}C_{0} & \dots & \prescript{1}{0}C_{k-1} & \dots \\
\vdots & \vdots & \ddots & \vdots & \dots \\
\prescript{0}{k-1}C & \prescript{1}{k-1}C_{0} & \dots & \prescript{1}{k-1}C_{k-1} & \dots \\
\prescript{0}{0}C \prescript{0}{0}C & \prescript{0}{0}C \prescript{1}{0}C_{0}+\prescript{0}{0}C \prescript{1}{0}C_{0} & \dots & \prescript{0}{0}C \prescript{1}{0}C_{k-1}+\prescript{0}{0}C \prescript{1}{0}C_{k-1} & \dots\\
\vdots & \vdots & \ddots & \vdots & \dots\\
\prescript{0}{0}C \prescript{0}{k-1}C & \prescript{0}{0}C \prescript{1}{k-1}C_{0}+\prescript{0}{k-1}C \prescript{1}{0}C_{0} & \dots & \prescript{0}{0}C \prescript{1}{k-1}C_{k-1}+\prescript{0}{k-1}C \prescript{1}{0}C_{k-1} & \dots \\
\prescript{0}{1}C \prescript{0}{0}C & \prescript{0}{1}C \prescript{1}{0}C_{0}+\prescript{0}{0}C \prescript{1}{1}C_{0} & \dots & \prescript{0}{1}C \prescript{1}{0}C_{k-1}+\prescript{0}{0}C \prescript{1}{1}C_{k-1} & \dots \\
\vdots & \vdots & \ddots & \vdots & \dots \\
\prescript{0}{1}C \prescript{0}{k-1}C & \prescript{0}{1}C \prescript{1}{k-1}C_{0}+\prescript{0}{k-1}C \prescript{1}{1}C_{0} & \dots & \prescript{0}{1}C \prescript{1}{k-1}C_{k-1}+\prescript{0}{k-1}C \prescript{1}{1}C_{k-1} & \dots\\
\vdots & \vdots & \ddots & \vdots & \dots \\
\prescript{0}{k-1}C \prescript{0}{k-1}C & \prescript{0}{k-1}C \prescript{1}{k-1}C_{0}+\prescript{0}{k-1}C \prescript{1}{k-1}C_{0} & \dots & \prescript{0}{k-1}C \prescript{1}{k-1}C_{k-1}+\prescript{0}{k-1}C \prescript{1}{k-1}C_{k-1} & \dots\\
\vdots & \vdots & \vdots & \vdots & \ddots\\
\end{array} \right]
\endgroup
\end{equation}
By replacing the groups of elements from (3.11) with the block-matrix notation of Kronecker powers of arrays ${\prescript{j}{p}C_{l_{1},...,l_{j}}}$, we arrive at (3.6).
\end{proof}
For ${k=1}$, the expression (3.6) reduces to (2.2). Furthermore, the analogue of the Theorem 2.1 can be proven. If ${\prescript{0}{p}C=0}$ for all ${p}$, then all of the Kronecker products involving ${\prescript{0}{p}C}$ vanish. Since all of the block matrices below the diagonal in (3.6) include the ${\prescript{0}{p}C}$ in the Kronecker product, the lower part of the block transition matrix vanishes and the matrix becomes block upper triangular. The block matrices on the diagonal of ${T_{ab}}$ are the Kronecker powers of the array standing near the linear term ${\prescript{1}{p}C_{l}}$. The Kronecker powers of ${\prescript{1}{p}C_{l}}$ are defined as in \cite{Book_8}:
\begin{equation}
\begingroup 
\setlength\arraycolsep{1.75pt}
\left(\prescript{1}{}{\bf C}\right)^{\alpha}=
\begin{bmatrix}
\prescript{1}{0}C_{0}\left(\prescript{1}{}{\bf C}\right)^{\alpha-1} & \prescript{1}{0}C_{1}\left(\prescript{1}{}{\bf C}\right)^{\alpha-1}  & \dots  & \prescript{1}{0}C_{k-1}\left(\prescript{1}{}{\bf C}\right)^{\alpha-1} \\
\prescript{1}{1}C_{0}\left(\prescript{1}{}{\bf C}\right)^{\alpha-1}  & \prescript{1}{1}C_{1}\left(\prescript{1}{}{\bf C}\right)^{\alpha-1}  & \dots  & \prescript{1}{1}C_{k-1}\left(\prescript{1}{}{\bf C}\right)^{\alpha-1} \\
\vdots & \vdots & \ddots & \vdots \\
\prescript{1}{k-1}C_{0}\left(\prescript{1}{}{\bf C}\right)^{\alpha-1}  & \prescript{1}{k-1}C_{1}\left(\prescript{1}{}{\bf C}\right)^{\alpha-1}  & \dots  & \prescript{1}{k-1}C_{k-1}\left(\prescript{1}{}{\bf C}\right)^{\alpha-1} \\
\end{bmatrix}
,\left(\prescript{1}{}{\bf C}\right)^{0}=1
\endgroup
\end{equation}
Therefore, the transition matrix can be brought to the upper triangular form if ${\forall p \ \prescript{0}{p}C=0 \wedge \forall p>l \ \prescript{1}{p}C_{l}=0}$. Even though in general recurrences do not fulfill those conditions, variables can be (in some cases) transformed in such a way that the new transition matrix is upper triangular.
\begin{lemma}
Let the recurrence on ${u_{i}^{k}}$ be defined as in (3.2). There exists a linear transformation of variables for ${F_{k}}$ such that the transition matrix ${T_{ab}}$ corresponding to the new sequence is upper triangular if and only if there exists ${k}$ dimensional vector ${B}$ satisfying:
\begin{equation}
B_{p}=\sum_{s=0}^{m}\frac{1}{s!}\sum_{l_{1},...,l_{s}=0}^{k-1}\sum_{\sigma}\prescript{s}{p}C_{l_{\sigma_{1}},...,l_{\sigma_{s}}}\prod_{t=1}^{s}B_{l_{t}}
\end{equation}
\end{lemma}
\begin{proof}
The general (invertible) linear transformation of ${k}$ variables ${u_{i}^{k}}$ can be written in the matrix form:
\begin{equation}
{\bf z}^{i}=A{\bf z}'^{i}+B\iff {\bf z}'^{i}=A^{-1}({\bf z}^{i}-B)
\end{equation}
with ${det(A)\neq 0}$. The new recurrence defined on primed variables ${{\bf z}'^{i}}$ satisfies:
\begin{equation}
z'^{i}_{p-1}=\prescript{0}{p}C'+\sum_{j=1}^{m}\sum_{l_{1},...,l_{j}=0}^{k-1}\prescript{j}{p}C'_{l_{1},...,l_{j}}\prod_{s=1}^{j}z'^{i-1}_{l_{s}}
\end{equation}
where:
\begin{equation}
\begin{aligned}
&\prescript{0}{k}C'=-A_{kj}^{-1}B_{j}+\sum_{s=0}^{m}\sum_{j,l_{1},...,l_{s}=0}^{k-1}A_{kj}^{-1}\prescript{s}{j}{\hat C}_{l_{1},...,l_{s}}\prod_{t=1}^{s}B_{l_{t}}\\
&\prescript{p}{k}C'_{x_{1},...,x_{p}}=\sum_{s=p}^{m}\binom{s}{p}\sum_{j,l_{1},...,l_{s}=0}^{k-1}A_{kj}^{-1}\prescript{s}{j}{\hat C}_{l_{1},...,l_{s}}\left(\prod_{t_{1}=1}^{p}A_{l_{t_{1}}x_{t_{1}}}\right)\left(\prod_{t_{2}=p+1}^{s}B_{l_{t_{2}}}\right)\\
&\prescript{s}{j}{\hat C}_{l_{1},...,l_{s}}=\frac{1}{s!}\sum_{\sigma}\prescript{s}{j}C_{l_{\sigma_{1}},...,l_{\sigma_{s}}}
\end{aligned}
\end{equation}
The symbol ${\prescript{s}{j}{\hat C}_{l_{1},...,l_{s}}}$ denotes the symmetrized array that is equal to ${\prescript{s}{j}C_{l_{1},...,l_{s}}}$ averaged over all of the permutations ${\sigma}$ of the right subscripts.\\\indent
It is trivial that convolution of linear transformations is a linear transformation. We firstly perform a shift of the variables by the vector ${B}$. From (3.16) we have:
\begin{equation}
\prescript{0}{p}C'=-B_{k}+\sum_{s=0}^{m}\frac{1}{s!}\sum_{l_{1},...,l_{s}=0}^{k-1}\sum_{\sigma}\prescript{s}{p}C_{l_{\sigma_{1}},...,l_{\sigma_{s}}}\prod_{t=1}^{s}B_{l_{t}}
\end{equation}
Setting ${\prescript{0}{p}C'=0}$, above equation reduces to a set of coupled polynomial equations. Unfortunately, the resulting system of equations is not always solvable. Assuming that (3.17) has a solution, the other arrays transform according to (3.16).\\\indent
Subsequently a linear transformation is performed with ${B=0}$. The constant term ${\prescript{0}{p}C''=\prescript{0}{p}C'=0}$ (i.e. still vanishes, regardless of the choice of ${A}$), while the linear term is:
\begin{equation}
\prescript{1}{p}C''_{t}=\sum_{j,l_{1}=0}^{k-1}A_{pj}^{-1}\prescript{1}{j}C'_{l}A_{lt}
\end{equation}
Similarly, there always exists a suitable choice of  ${k\times k}$ dimensional matrix ${A}$ such that ${\prescript{1}{p}C_{t}}$ becomes upper triangular in lower indices \cite{Book_8} and thus the transition matrix can be brought to an upper triangular form provided (3.13) has a solution.
\end{proof}
If the shift vector ${B}$ satisfying (3.13) exists, its not necessary unique. The situation is analogous to the one found for a uni-variable recurrence. Different choices of ${B}$ satisfying (3.13) with ${\prescript{0}{p}C'=0}$ will lead to ${\prescript{1}{p}C_{l}}$ with different eigenvalues. Therefore, it is suitable to choose ${B}$ such that it satisfies the diagonalizability conditions stated in Lemma 3.3.\\\indent
It is worth noting that even though the arrays ${\prescript{j}{p}C_{l_{1},...,l_{j}}}$ vanish for a non-decreasing order of lower right indices, the same may not be true for the transformed arrays ${\prescript{j}{p}C'_{l_{1},...,l_{j}}}$. Therefore, after transforming the coefficient  arrays with the lineary transformation (3.13), it is necessary to upper triangularize the transformed coefficients so that they also satisfy ${\prescript{j}{p}C'_{l_{1},...,l_{j}}=0 \ if \ \neg (l_{1}\leq ... \leq l_{j})}$.\\\indent
Since the eigenvalues of the Kronecker product of two matrices are the products of eigenvalues of the two matrices, we have:
\begin{equation}
Spec(\left(\prescript{1}{}{\bf C}\right)^{\alpha})=Spec(\prescript{1}{}{\bf C})^{\alpha}
\end{equation}
where the right-hand side represents the Cartesian powers of the set of eigenvalues of ${\prescript{1}{}{\bf C}}$ \cite{Book_11}. As a result, the general eigenvalue is of the form:
\begin{equation}
\lambda=\lambda'^{a_{0}}_{0}\lambda'^{a_{1}}_{1}...\lambda'^{a_{k-1}}_{k-1}
\end{equation}
for some numbers ${a_{p}\in \mathbb{N}, \lambda_{p}'\in Spec(\prescript{1}{}{\bf C})}$. In the same way as in (3.7), the l-th eigenvalue of the transition can be written as:
\begin{equation}
\lambda_{l}=
\begin{cases}
1 \ for \ l=0 \\
\prod_{s=1}^{\lfloor \log_{k}(1+l(k-1)) \rfloor}\lambda'_{{\lfloor (l-\frac{1-k^{s}}{1-k})k^{1-s}\rfloor\%k}}
\end{cases}
\end{equation}
The matrix is diagonalizable if all of its eigenvalues are different, thus a sufficient condition for the existence of diagonal form can be constructed.\\\indent
Although this is true, a problem arises due to redundancy in description in terms of the auxiliary vectors (3.4). Since ${{\bf y}^{i}}$ are composed of Kronecker powers of ${{\bf z}^{i}}$, there are components that are redundant, i.e. that represent the same product of variables ${u_{i}^{k}}$. For example, if ${k=2}$, the infinite auxiliary vectors are:
\begin{equation}
{\bf y}^{i}=[1,u_{i}^{1},u_{i}^{2},(u_{i}^{1})^{2},u_{i}^{1}u_{i}^{2},u_{i}^{2}u_{i}^{1},(u_{i}^{2})^2,...]^{T}
\end{equation}
and it is apparent that the 4-th and 5-th components are the same.\\\indent
Due to the redundancy, the repeating copies of the same eigenvalues appear in the upper triangular form of the transition matrix (3.6). Lemmas 2.4 and 2.5 are applicable only if the eigenvalues of ${T_{ab}}$ are non-degenerate, hence the Lemmas cannot be used to diagonalize the transition matrix. Nevertheless, the special choice of ${\prescript{j}{p}C_{l_{1},...,l_{j}}}$ allows one to "reduce" the transition matrix and get rid of the undesired copies.\\\indent
The redundant terms arise from the permutation of lower right indices of ${\prescript{j}{p}C'_{l_{1},...,l_{j}}}$. By definition (the begining of Section 3) the arrays vanish for a decreasing sequence of lower indices. All of the columns in ${T_{ab}}$ associated with the copies of the same power of variables ${u_{i}^{k}}$ vanish except for the diagonal element. For example, since the 4-th and 5-th element of the auxiliary vectors (3.22) are similar, all of the elements ${T_{a5}=0}$, ${a\neq 5}$ vanish. This observation allows one to get rid of the redundant columns (and rows) of ${T_{ab}}$, hence making the calculation of the powers of the transition matrix feasible. In other words, the redundancy has to be taken into account when assessing the diagonalizability of ${T_{ab}}$.\\\indent
In order to differentiate between quantities related to the redundant (original) and non-redundant system, tilde notation is introduced. The quantities related to the system without redundancy (vectors without redundant components and matrices without redundant rows and columns) are denoted with tilde above the symbol. For example, the equivalent of the auxiliary vector (3.22) for the non-redundant Carleman embedding is:
\begin{equation}
{\bf \tilde{y}}^{i}=[1,u_{i}^{1},u_{i}^{2},(u_{i}^{1})^{2},u_{i}^{1}u_{i}^{2},(u_{i}^{2})^2,...]^{T}
\end{equation}
Similarly, the transition matrix for ${{\bf \tilde{y}}^{i}}$ is denoted by ${\tilde{T}_{ab}}$ respectively. Although it is the non-redundant Carleman embedding defined by ${\tilde{T}_{ab}}$ and ${{\bf \tilde{y}}^{i}}$ that is eventually diagonalized, it is hard to express it using a concise formula. Equivalents of the expression (3.4) for ${{\bf \tilde{y}}^{i}}$ and (3.6) for ${\tilde{T}_{ab}}$ has not been found. As opposed to the redundant system, the transition matrix ${\tilde{T}_{ab}}$ cannot be simply expressed in terms of Kronecker products. On the other hand, it is much easier to calculate powers of ${\tilde{T}_{ab}}$ rather than ${T_{ab}}$. The following can be said about diagonalizability of the "non-redundant" transition matrix:
\begin{lemma}
A sufficient condition for the infinite upper triangular transition matrix ${\tilde{T}_{ab}}$ to be diagonalizable is for the corresponding ${k\times k}$ dimensional matrix ${\prescript{1}{}{\bf C}}$ to have eigenvalues ${\lambda_{p}'}$ such that the equation:
\begin{equation}
\begin{cases}
x_{0}+x_{1}\log_{|\lambda_{0}|}{|\lambda_{1}|}+...+x_{k-1}\log_{|\lambda_{0}|}{|\lambda_{k-1}|}=0\\
Arg(\lambda_{0})x_{0}+Arg(\lambda_{1})x_{1}+...+Arg(\lambda_{k-1})x_{k-1}=2\pi t
\end{cases}
\end{equation}
has no (non-trivial) solution ${x_{0},...,x_{k-1},t\in \mathbb{Z}}$. In case of two variables ${k=2}$, this can be further simplified to:
\begin{equation}
    \log_{|\lambda_{0}|}{|\lambda_{1}|}\not \in \mathbb{Q} \lor \frac{\pi}{Arg(\lambda_{1})-Arg(\lambda_{0})\log_{|\lambda_{0}|}{|\lambda_{1}|}}\not \in \mathbb{Q}
\end{equation}
\end{lemma}
\begin{proof}
As stated before, the redundancy in the auxiliary vectors ${{\bf y}^{i}}$ implies that there always will be an eigenvector with algebraic multiplicity of at least two. However, since all of the elements of the corresponding column except the diagonal vanish, a sufficient condition for diagonalization can be still formulated based on the tilded transition matrix ${\tilde{T}_{ab}}$.\\\indent
The condition for ${\tilde{T}_{ab}}$ to have different eigenvalues can be written as:
\begin{equation}
\lambda'^{a_{0}}_{0}\lambda'^{a_{1}}_{1}...\lambda'^{a_{k-1}}_{k-1}=\lambda'^{a_{0}'}_{0}\lambda'^{a_{1}'}_{1}...\lambda'^{a_{k-1}'}_{k-1}\implies \forall p \ a_{p}=a_{p}'
\end{equation}
Dividing by the primed side:
\begin{equation}
\lambda'^{a_{0}-a_{0}'}_{0}\lambda'^{a_{1}-a_{1}'}_{1}...\lambda'^{a_{k-1}-a_{k-1}'}_{k-1}=1
\end{equation}
Therefore, the condition is satisfied if and only if the only solution ${a_{p}-a_{p}'\in \mathbb{Z}}$ to (3.27) is ${a_{p}-a_{p}'=0}$ for all ${p}$. By rewriting the eigenvalues in the polar form ${\lambda_{p}'=r_{p}e^{i\theta_{p}}}$:
\begin{equation}
\begin{cases}
r_{0}^{a_{0}-a_{0}'}r_{1}^{a_{1}-a_{1}'}...r_{k-1}^{a_{k-1}-a_{k-1}'}=1\\
e^{i\theta_{0}(a_{0}-a_{0}')}e^{i\theta_{1}(a_{1}-a_{1}')}...e^{i\theta_{k-1}(a_{k-1}-a_{k-1}')}=1
\end{cases}
\end{equation}
and taking the logarithm of both sides of (3.28), we obtain a set of linear equations:
\begin{equation}
\begin{cases}
(a_{0}-a_{0}')+(a_{1}-a_{1}')\log_{r_{0}}{r_{1}}+...+(a_{k-1}-a_{k-1}')\log_{r_{0}}{r_{k-1}}=0\\
\theta_{0}(a_{0}-a_{0}')+\theta_{1}(a_{1}-a_{1}')+...+\theta_{k-1}(a_{k-1}-a_{k-1}')=2\pi t, t\in \mathbb{Z}
\end{cases}
\end{equation}
There always exists a trivial solution ${\forall p \ a_{p}-a_{p}'=0}$.\\\indent
In case ${k=2}$, the equations reduce to:
\begin{equation}
\begin{cases}
(a_{0}-a_{0}')+(a_{1}-a_{1}')\log_{r_{0}}{r_{1}}=0\\
\theta_{0}(a_{0}-a_{0}')+\theta_{1}(a_{1}-a_{1}')=2\pi t, t\in \mathbb{Z}
\end{cases}
\end{equation}
The upper equation has integer solutions if and only if ${\log_{r_{0}}{r_{1}}\in \mathbb{Q}}$. Substituting the result to the lower equation, we conclude that the above system of equations does not have an integer solution if:
\begin{equation}
\log_{r_{0}}{r_{1}}\not \in \mathbb{Q} \lor \frac{\pi}{\theta_{1}-\theta_{0}\log_{r_{0}}{r_{1}}}\not \in \mathbb{Q}
\end{equation}
\end{proof}
If the transition matrix ${\tilde{T}_{ab}}$ can be brought into an upper triangular form, Lemmas 2.4 and 2.5 can be used to diagonalize it. This allows one to give an explicit expression for powers of ${\tilde{T}_{ab}}$, thus solving the recurrence. We now prove the main theorem of the article:

\begin{theorem}
Let ${u_{i}^{n}\in \mathbb{C}, i\geq 0, 1\leq n\leq k}$ be ${k}$ sequences of numbers satisfying the recurrence relation (3.2) and let ${A}$ and ${B}$ denote the parameters of a linear transformation that brings the transition matrix to an upper triangular form, i.e. ${k\times k}$ matrix ${A}$ and k dimensional vector ${B}$ such that:
\begin{equation}
B_{p}=\sum_{s=0}^{m}\frac{1}{s!}\sum_{l_{1},...,l_{s}=0}^{k-1}\sum_{\sigma}\prescript{s}{p}C_{l_{\sigma_{1}},...,l_{\sigma_{s}}}\prod_{t=1}^{s}B_{l_{t}}
\end{equation}
Furthermore, let the ${k\times k}$ array ${\prescript{1}{p}C''_{x}}$ be upper triangular, where:
\begin{equation}
\begin{aligned}
&\prescript{1}{p}C''_{x}=\sum_{j,l_{1}=0}^{k-1}A_{pj}^{-1}\prescript{1}{j}C'_{l}A_{lx}\\
&\prescript{1}{p}C'_{l_{1}}=\sum_{s=1}^{m}\sum_{l_{2},...,l_{s}=0}^{k-1}\frac{1}{(s-1)!}\sum_{\sigma}\prescript{s}{p}C_{l_{\sigma_{1}},...,l_{\sigma_{s}}}\left(\prod_{t=2}^{s}B_{l_{t}}\right)\\
\end{aligned}
\end{equation}
Let the eigenvalues ${\prescript{1}{p}C''_{x}}$ (by eigenvalues of the array we understand the eigenvalues of the corresponding matrix created by considering only two lower indices) ${\lambda_{p}'}$ be such that the equation:
\begin{equation}
\begin{cases}
x_{0}+x_{1}\log_{|\lambda_{0}|}{|\lambda_{1}|}+...+x_{k-1}\log_{|\lambda_{0}|}{|\lambda_{k-1}|}=0\\
Arg(\lambda_{0})x_{0}+Arg(\lambda_{1})x_{1}+...+Arg(\lambda_{k-1})x_{k-1}=2\pi t
\end{cases}
\end{equation}
has no non-trivial solutions ${x_{0},...,x_{k-1},t\in \mathbb{Z}}$. Then the solution to the recurrence is given in terms of the initial condition ${{\bf \tilde{y}}^{0}}$ by:
\begin{equation}
u_{i}^{k}=B_{k-1}+\sum_{s=1}^{k}\sum_{l=1}^{\infty}\sum_{j=1}^{l}A_{k-1,s-1}\tilde{v}'^{s}_{j}\tilde{T}'^{i}_{jj}\tilde{P}'^{-1}_{jl}\tilde{y}'^{0}_{l}
\end{equation}
where:
\begin{equation}
\begin{aligned}
&P'^{-1}_{l_{0},l_{p+1}}=
\begin{cases}
0 \ for \ l_{0}>l_{p+1} \\
1 \ for \ l_{0}=l_{p+1} \\
\sum_{p=0}^{l_{p+1}-l_{0}-1} \ \sum_{l_{0}<l_{1}<...<l_{p}<l_{p+1}}(-1)^{p+1}\prod_{j=0}^{p}\frac{\tilde{v}'^{l_{j}}_{l_{j+1}}}{\tilde{v}'^{l_{j}}_{l_{j}}}\  
for \ l_{0}<l_{p+1}
\end{cases}
\\
&v'^{l_{0}}_{l_{p+1}}=
\begin{cases}
0 \ for \ l_{0}>l_{p+1} \\
1 \ for \ l_{0}=l_{p+1} \\
\sum_{p=0}^{l_{p+1}-l_{0}-1} \ \sum_{l_{0}<l_{1}<...<l_{p}<l_{p+1}}(-1)^{p+1}\prod_{j=0}^{p}\frac{\tilde{T}'_{l_{j},l_{j+1}}}{\tilde{T}'_{l_{j},l_{j}}-\lambda'_{l_{p+1}}} \ for \ l_{0}<l_{p+1}

\end{cases}
\end{aligned}
\end{equation}
\end{theorem}
\begin{proof}
Conditions (3.32) and (3.33) ensure that after the linear transformation ${{\bf z}'^{i}=A^{-1}({\bf z}^{i}-B)}$ the corresponding transition matrix ${T_{ab}'}$ is upper triangular. Similarly, the condition (3.34) is simply a sufficient condition for the diagonalization of ${\tilde{T'}_{ab}}$ stated in Lemma 3.3.\\\indent
Starting from Theorem 3.1, we have ${u'^{k}_{i}={{\bf e}_{k}}T'^{i}{{\bf y}^{0}}}$. However, as mentioned before, due to redundancy of description in terms of ${T_{ab}}$, the corresponding ${\tilde{T}_{ab}}$ is used instead. Furthermore, since there is no redundancy in the first ${k+1}$ elements of ${\tilde{\bf y}^{i}}$, the solution can be equivalently written as:
\begin{equation}
u'^{k}_{i}={{\bf e}_{k}}T'^{i}{{\bf y}^{0}}={{\bf e}_{k}}\tilde{T}'^{i}{{\bf \tilde{y}}^{0}}=\sum_{l=1}^{\infty}\sum_{j=1}^{l}\tilde{v}'^{k}_{j}\tilde{T}'^{i}_{jj}\tilde{P}'^{-1}_{jl}\tilde{y}'^{0}_{l}
\end{equation}
Where the symbols with tilde are symbols related to ${\tilde{T}_{ab}}$ free from the redundancy mentioned before. The prime refers to the quantities related to the recurrence obtained after the linear transformation of variables. Transforming (3.37) back to the original sequences:
\begin{equation}
u_{i}^{k}=\sum_{s=1}^{k}A_{k-1,s-1}u'^{s}_{i}+B_{k-1}=B_{k-1}+\sum_{s=1}^{k}\sum_{l=1}^{\infty}\sum_{j=1}^{l}A_{k-1,s-1}\tilde{v}'^{s}_{j}\tilde{T}'^{i}_{jj}\tilde{P}'^{-1}_{jl}\tilde{y}'^{0}_{l}
\end{equation}
Thus, we arrive at the expression (3.35). The eigenvectors and the inverse of the infinite modal matrix ${P^{-1}}$ are given in Lemmas 2.4 and 2.5.
\end{proof}
The solutions to the recurrence (3.2) have a general form:
\begin{equation}
u_{i}=\sum_{j_{1},...,j_{k}=0}^{\infty}f_{j_{1},...,j_{k}}(i)\prod_{s=1}^{k}\left(u_{0}^{s}\right)^{j_{s}}
\end{equation}
which closely resembles the formula (2.37). A noticeable difference is that for the multi-variable recurrence, the solution takes the form of a multi-variable (instead of uni-variable) power series in the initial condition. 
\begin{example}
As an example, we consider the simple recurrence defined by:
\begin{equation}
\begin{cases}
u_{i}^{1}=8u_{i-1}^{1}+10u_{i-1}^{2}+(u_{i-1}^{1})^{2}+3u_{i-1}^{1}u_{i-1}^{2}+(u_{i-1}^{2})^{2}\\
u_{i}^{2}=-3u_{i-1}^{1}-3u_{i-1}^{2}+(u_{i-1}^{1})^{2}-u_{i-1}^{1}u_{i-1}^{2}+(u_{i-1}^{2})^{2}
\end{cases}
\end{equation}
The arrays ${\prescript{j}{p}C_{l_{1},...,l_{j}}}$ can be obtained by comparing (3.40) with (3.2):
\begin{equation}
\prescript{0}{p}C=
\begin{bmatrix}
0 \\
0 \\
\end{bmatrix}
,\prescript{1}{p}C_{l}=
\begin{bmatrix}
8 & 10 \\
-3 & -3 \\
\end{bmatrix}
,\prescript{2}{0}C_{l_{1},l_{2}}=
\begin{bmatrix}
1 & 3 \\
0 & 1 \\
\end{bmatrix}
,\prescript{2}{1}C_{l_{1},l_{2}}=
\begin{bmatrix}
1 & -1 \\
0 & 1 \\
\end{bmatrix}
\end{equation}
In order to bring the transition matrix into the upper triangular form, the variables are transformed by a ${2\times 2}$ matrix ${A}$:
\begin{equation}
A=
\begin{bmatrix}
1 & 2 \\
-3 & -  5 \\
\end{bmatrix}
\end{equation}
to obtain the new recurrence:
\begin{equation}
\begin{cases}
u'^{1}_{i}=2u'^{1}_{i-1}+87(u'^{1}_{i-1})^{2}+67u'^{1}_{i-1}u'^{2}_{i-1}+13(u'^{2}_{i-1})^{2}\\
u'^{2}_{i}=3u'^{2}_{i-1}-212(u'^{1}_{i-1})^{2}-164u'^{1}_{i-1}u'^{2}_{i-1}-32(u'^{2}_{i-1})^{2}
\end{cases}
\end{equation}
with ${\prescript{0}{p}C'=0}$ and:
\begin{equation}
\prescript{1}{p}C'_{l}=
\begin{bmatrix}
2 & 0 \\
0 & 3 \\
\end{bmatrix},
\prescript{2}{0}C'_{l_{1},l_{2}}=
\begin{bmatrix}
87 & 67 \\
0 & 13 \\
\end{bmatrix},
\prescript{2}{1}C'_{l_{1},l_{2}}=
\begin{bmatrix}
-212 & -164 \\
0 & -32 \\
\end{bmatrix}
\end{equation}
The eigenvalues satisfy the conditions (3.34), therefore ${\tilde{T}'_{ab}}$ can be diagonalized and the powers of the transition matrix can be calculated. The infinite transition matrix is:
\begin{equation}
T_{ab}=
\left[ \begin{array}{@{}*{20}{c}@{}}
1 & 0 & 0 & 0 & 0 & 0 & 0 & 0 & 0 & 0 & 0 & 0 & \dots \\
0 & 2 & 0 & 87 & 67 & 0 & 13 & 0 & 0 & 0 & 0 & 0 & \dots \\
0 & 0 & 3 & -212 & -164 & 0 & -32 & 0 & 0 & 0 & 0 & 0 & \dots \\
0 & 0 & 0 & 4 & 0 & 0 & 0 & 348 & 268 & 0 & 52 & 0 & \dots \\
0 & 0 & 0 & 0 & 6 & 0 & 0 & -424 & 67 & 0 & 137 & 0 & \dots \\
0 & 0 & 0 & 0 & 0 & 6 & 0 & -424 & 67 & 0 & 137 & 0 & \dots \\
0 & 0 & 0 & 0 & 0 & 0 & 9 & 0 & -1272 & 0 & -987 & 0 & \dots \\
0 & 0 & 0 & 0 & 0 & 0 & 0 & 8 & 0 & 0 & 0 & 0 & \dots \\
0 & 0 & 0 & 0 & 0 & 0 & 0 & 0 & 12 & 0 & 0 & 0 & \dots \\
0 & 0 & 0 & 0 & 0 & 0 & 0 & 0 & 0 & 12 & 0 & 0 & \dots \\
0 & 0 & 0 & 0 & 0 & 0 & 0 & 0 & 0 & 0 & 18 & 0 & \dots \\
0 & 0 & 0 & 0 & 0 & 0 & 0 & 0 & 0 & 0 & 0 & 12 & \dots \\
\vdots & \vdots & \vdots & \vdots & \vdots & \vdots & \vdots & \vdots & \vdots & \vdots & \vdots & \vdots & \ddots \\
\end{array} \right]
\end{equation}
The 5-th, 9-th and 11-th rows and columns are examples of the arising redundancy. The Lemmas 2.4 and 2.5 can now be applied to the reduced transition matrix ${\tilde{T}_{ab}}$, which is obtained by eliminating all of the redundant rows and columns:
\begin{equation}
\tilde{T}_{ab}=\\
\begin{bmatrix}
1 & 0 & 0 & 0 & 0 & 0 & \dots \\
0 & 2 & 0 & 87 & 67 & 13 & \dots \\
0 & 0 & 3 & -212 & -164 & -32 & \dots \\
0 & 0 & 0 & 4 & 0 & 0 & \dots \\
0 & 0 & 0 & 0 & 6 & 0 & \dots \\
0 & 0 & 0 & 0 & 0 & 9 & \dots \\
\vdots & \vdots & \vdots & \vdots & \vdots & \vdots & \ddots \\
\end{bmatrix}
\end{equation}
This can be further rewritten in the transformed form:

\begin{equation}
\begin{aligned}
\tilde{T}_{ab}&=
\begin{bmatrix}
1 & 0 & 0 & 0 & 0 & 0 & \dots \\
0 & 1 & 0 & 87 & 201 & 39 & \dots \\
0 & 0 & 1 & -424 & -656 & -112 & \dots \\
0 & 0 & 0 & 2 & 0 & 0 & \dots \\
0 & 0 & 0 & 0 & 12 & 0 & \dots \\
0 & 0 & 0 & 0 & 0 & 21 & \dots \\
\vdots & \vdots & \vdots & \vdots & \vdots & \vdots & \ddots \\
\end{bmatrix}
\begin{bmatrix}
1 & 0 & 0 & 0 & 0 & 0 & \dots \\
0 & 2 & 0 & 0 & 0 & 0 & \dots \\
0 & 0 & 3 & 0 & 0 & 0 & \dots \\
0 & 0 & 0 & 4 & 0 & 0 & \dots \\
0 & 0 & 0 & 0 & 6 & 0 & \dots \\
0 & 0 & 0 & 0 & 0 & 9 & \dots \\
\vdots & \vdots & \vdots & \vdots & \vdots & \vdots & \ddots \\
\end{bmatrix}
\times \\
&\times
\begin{bmatrix}
1 & 0 & 0 & 0 & 0 & 0 & \dots \\
0 & 1 & 0 & -\frac{87}{2} & -\frac{67}{4} & -\frac{13}{7} & \dots \\[2pt]
0 & 0 & 1 & 212 & \frac{164}{3} & \frac{16}{3} & \dots \\
0 & 0 & 0 & \frac{1}{2} & 0 & 0 & \dots \\
0 & 0 & 0 & 0 & \frac{1}{12} & 0 & \dots \\
0 & 0 & 0 & 0 & 0 & \frac{1}{21} & \dots \\
\vdots & \vdots & \vdots & \vdots & \vdots & \vdots & \ddots \\
\end{bmatrix}
\end{aligned}
\end{equation}
The matrices to the left and right of the diagonal matrix are ${P}$ and ${P^{-1}}$ respectively. The powers of ${\tilde{T}_{ab}}$ can now be easily calculated. Using (3.5), we finally arrive at:
\begin{equation}
\begin{aligned}
&u'^{1}_{i}=2^{i}u'^{1}_{0}+\frac{87}{2}(4^{i}-2^{i})\left(u'^{1}_{0}\right)^{2}+\frac{67}{4}(6^{i}-2^{i})u'^{1}_{0}u'^{2}_{0}+\frac{13}{7}(9^{i}-2^{i})\left(u'^{2}_{0}\right)^{2}+...\\
&u'^{2}_{i}=3^{i}u'^{2}_{0}-212(4^{i}-3^{i})\left(u'^{1}_{0}\right)^{2}-\frac{164}{3}(6^{i}-3^{i})u'^{1}_{0}u'^{2}_{0}-\frac{16}{3}(9^{i}-3^{i})\left(u'^{2}_{0}\right)^{2}+...
\end{aligned}
\end{equation}
The above solution to the primed recurrence (3.43) can now be transformed back to the original variables, therefore giving a rather complicated expression for the solution to the recurrence (3.40):
\begin{equation}
\begin{aligned}
u_{i}^{1}&=\left(6\cdot 3^{i}-5\cdot 2^{i}\right)u_{0}^{1}+10\left(3^{i}-2^{i}\right)u_{0}^{2}+\\
&+\left(\frac{87}{7}9^{i}-\frac{307}{4}6^{i}+\frac{413}{2}4^{i}-192\cdot 3^{i}+\frac{1395}{28}2^{i}\right)\left(u_{0}^{1}\right)^{2}+
\\&+\left(\frac{290}{7}9^{i}-\frac{3377}{12}6^{i}+826\cdot 4^{i}-\frac{2440}{3}3^{i}+\frac{6365}{28}2^{i}\right)u_{0}^{1}u_{0}^{2}+
\\&+\left(\frac{725}{21}9^{i}-\frac{1535}{6}6^{i}+826\cdot 4^{i}-\frac{2608}{3}3^{i}+\frac{3705}{14}2^{i}\right)\left(u_{0}^{2}\right)^{2}+...\\
u_{i}^{2}&=-3\left(3^{i}-2^{i}\right)u_{0}^{1}-\left(5\cdot 3^{i}-6\cdot 2^{i}\right)u_{0}^{2}+\\
&+\left(\frac{15}{7}9^{i}+\frac{53}{4}6^{i}-\frac{163}{2}4^{i}+96\cdot 3^{i}-\frac{837}{28}2^{i}\right)\left(u_{0}^{1}\right)^{2}+\\
&+\left(\frac{50}{7}9^{i}+\frac{583}{12}6^{i}-326\cdot 4^{i}+\frac{1220}{3}3^{i}-\frac{3819}{28}2^{i}\right)u_{0}^{1}u_{0}^{2}+\\
&+\left(\frac{125}{21}9^{i}+\frac{265}{6}6^{i}-326\cdot 4^{i}+\frac{1304}{3}3^{i}-\frac{2223}{14}2^{i}\right)\left(u_{0}^{2}\right)^{2}+...
\end{aligned}
\end{equation}
Comparing (3.49) with the general form of the solution (3.39), the coefficient functions can be obtained.
\end{example}

\section{Notes on multi-variable recurrences of arbitrary depth}
The Theorem 3.4 is also applicable to arbitrary-depth recursions. For any finite-depth recursion, auxiliary variables can be introduced and the original recurrence can be reduced to multi-variable recurrence of depth one. The method follows closely the one used for the finite-depth linear recurrences \cite{Book_2}.\\\indent
Let the depth-n recurrence be in the form (1.1). Set of ${k\times n}$ auxiliary variables is defined by:
\begin{equation}
u_{i-j}^{l}=u_{i}^{l+jk}
\end{equation}
for ${0\leq j<n}$. The equation (1.1) now takes the form:
\begin{equation}
\begin{cases}
u_{i}^{1} & =F_{1}\left(u_{i-1}^{1},..., u_{i-1}^{kn}\right)\\
&\vdots\\
u_{i}^{k} & =F_{k}\left(u_{i-1}^{1},..., u_{i-1}^{kn}\right)\\
u_{i}^{k+1} & =u_{i}^{1} \\
&\vdots\\
u_{i}^{kn} & =u_{i}^{k(n-1)}\\
\end{cases}
\end{equation}
Which is a multi-variable depth one recurrence, already discussed in Section 3.
\section{Further remarks}
In Section 2, an explicit formula for the solution of wide range of uni-variable depth-one recurrences is given. Although the final expressions is rather complicated, it can be used to derive the coefficient functions ${f_{j}(i)}$ in (2.37). Even if the recurrence does not fulfill the diagonalizability condition, there is a chance that it can be reformulated in terms of new variables. The transformation should be chosen such that the corresponding transition matrix has different eigenvalues and Lemmas 2.4 and 2.5 can be applied. The shift transformation is of particular importance, nevertheless any transformation that leaves the recurrence polynomial could be used. Transformations by a higher degree invertible polynomial, such us ${u_{i}'=u_{i}^3+1}$ are just one of many examples.\\\indent
Although the method allows to tackle many different polynomial recurrence relations, there are some notable shortcomings. The diagonalizability conditions in Sections 2 and 3 are only sufficient and not necessary for the transition matrix to be diagonalizable. Therefore, the method may not be applicable for all systems that are solvable by diagonalization approach. In addition, transition matrices which are known to be non-diagonalizabe are yet an another group of cases where Theorems 2.6 and 3.4 cannot be used. Possibly, an analogous method could be formulated, where Jordanization is used to calculate the powers of the transition matrix instead of diagonalization.\\\indent
The Theorem 3.4 utilizes the non-redundant system (with ${\tilde{T}}$) instead of the one defined in Theorem 3.1. It is dictated by the fact that the transition matrix defined for the auxiliary vectors ${3.4}$ necessary contains redundant copies of its eigenvalues. This may lead to difficulties with evaluating the expansion (3.39) for large ${j_{1},...,{j_{k}}}$, as only the general rule for constructing ${\tilde{T}}$ is given as opposed to explicit expression as in case of ${T}$, i.e. (3.9).\\\indent
It can be noted that examples 2.7 and 3.5 belong to a family of recurrences for which no shift is needed. The logistic map automatically fulfills the assumptions for ${r\neq 1}$, while the recurrence relation (3.40) needs to be transformed only by matrix ${A}$, with ${B=0}$. Those examples were chosen on purpose, as the solution of the transformed series can be truncated and the calculations simplify significantly. In general however, if one wants to obtain the exact solution to the original recurrence with non-zero shift parameter, all of the terms of the solution corresponding to the transformed recurrence have to be kept, which often makes the calculations cumbersome.\\\indent
Even though the Carleman linearization method does not provide a way to systematically approach the topic of polynomial recurrences, a large group of relations (1.1) fulfill the conditions and the method can be applied. Furthermore, it is straightforward to generalize the method to systems of recurrences (Theorem 3.4), as well as systems of arbitrary-depth recurrences, as discussed in Section 4.\\\indent
Improvements to the Carleman linearization of polynomial recurrences may include research on the diagonalizability of the transition matrices, possibly resulting in less-restrictive diagonalizability conditions. In addition, polynomial transformations may play an important role in further extending the group of recurrences the method is applicable to.
\section{Acknowledgements}
The author is grateful to Daniel Gagliardi for valuable remarks, and to the anonymous referees for devoting their time and effort to the manuscript.

\bibliographystyle{ieeetr}
\bibliography{References}

\end{document}